\documentclass[hyper,letterpaper]{JHEP3}
\usepackage{amsmath,amssymb,amsfonts,bm}
\usepackage{cite}
\usepackage{graphicx}
\usepackage{multirow}
\usepackage{verbatim}
\usepackage{appendix}

\usepackage{url}
\usepackage{float}

\newcommand{\bea}{\begin{eqnarray}}
\newcommand{\eea}{\end{eqnarray}}
\newcommand{\be}{\begin{equation}}
\newcommand{\ee}{\end{equation}}

\newcommand{\unknot}{{\,\raisebox{-.08cm}{\includegraphics[width=.37cm]{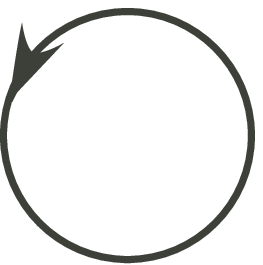}}\,}}

\newcommand{\R}{{\mathbb R}}
\newcommand{\C}{{\mathbb C}}

\newcommand{\cN}{{\mathcal{N}}}
\newcommand{\Li}{{\rm Li}}

\def\Tr{{\rm Tr \,}}

\def\k{\kappa}

\newcommand{\cM}{{\cal M }}
\newcommand{\cC}{{\cal C }}

\newcommand{\cO}{{\cal O }}
\newcommand{\cH}{{\cal H }}

\renewcommand{\P}{{\cal P}}
\newcommand{\cp}{{\mathbb{C}}{\mathbf{P}}}
\renewcommand{\S}{{\bf S}}

\renewcommand{\bar}{\overline}
\renewcommand{\hat}{\widehat}

\title{Super-A-polynomial}

\author{Hiroyuki Fuji$^{1,2}$  
and Piotr Su{\l}kowski$^{1,3}$
\\
$^1$ California Institute of Technology, Pasadena, CA 91125, USA \\
$^2$ Nagoya University, Dept. of Physics, Graduate School of Science, \\
$\ $ Furo-cho, Chikusa-ku, Nagoya 464-8602, Japan \\
$^3$ Faculty of Physics, University of Warsaw, ul. Ho{\.z}a 69, 00-681 Warsaw, Poland}

\abstract{We review a construction of a new class of algebraic curves, called super-$A$-polynomials, 
and their quantum generalizations. 
The super-$A$-polynomial is a two-parameter deformation of the $A$-polynomial known from knot theory
or Chern-Simons theory with $SL(2,\mathbb{C})$ gauge group.
The two parameters of the super-$A$-polynomial encode, respectively, the $t$-deformation which leads to the ``refined $A$-polynomial'',
and the $Q$-deformation which leads to the augmentation polynomial of knot contact homology.
For a given knot, the super-$A$-polynomial encodes the asymptotics of the corresponding $S^r$-colored HOMFLY homology for large $r$, while
the quantum super-$A$-polynomial provides recursion relations for such homology theories for each $r$.
The super-$A$-polynomial also admits a simple physical interpretation
as the defining equation for the space of SUSY vacua
in a circle compactification of the effective 3d $\cN=2$ theory associated to a given knot (complement).
We discuss properties of super-$A$-polynomials and illustrate them in many examples. 
\\
\\
\\
\\
\\
\\
{\tt CALT-68-2918}}

\begin{document}


\section{Introduction}
\label{sec:intro}

In past two decades remarkable relations between quantum field theory, string theory and knot theory have been found, following the seminal work by Witten \cite{Witten_Jones}. Among the others, such important mathematical developments as polynomial knot invariants, volume conjectures, $A$-polynomials, homological knot invariants, and more, have been interpreted from the perspective of high energy physics. In this note we summarize a construction of a new object in this line of research, the so-called super-$A$-polynomial, introduced in \cite{superA}. The super-$A$-polynomial can be regarded as a two-parameter generalization of an ordinary $A$-polynomial \cite{CCGL,Garoufalidis,Apol}. One of these parameters encodes information about the $t$-deformation of knot invariants arising upon categorification. The one-parameter generalization of an $A$-polynomial, depending only on $t$, has been introduced in \cite{FGS} as the so-called ``refined $A$-polynomial''. This parameter is related to the categorification of knot invariants and knot homologies, such as Khovanov homology \cite{Khovanov}, Khovanov-Rozansky homology \cite{KhR1}, or HOMFLY homology \cite{DGR}; more precisely, in this paper $t$ appears upon taking the Poincar\'e characteristis of the (conjectural) colored HOMFLY homology \cite{GS}. The second parameter, denoted by $a$ in this note, is related to Chern-Simons theory with $SU(N)$ gauge group. It also corresponds to the so-called $Q$-deformation of the $A$-polynomial introduced in \cite{AVqdef}, as well as the augmentation polynomial of knot contact homology \cite{NgFramed}. The super-$A$-polynomial captures information about both $a$ and $t$ at once, and, among the others, 
it encodes the asymptotics of the corresponding $S^r$-colored HOMFLY homologies for large $r$. In addition, its quantum deformation, the so-called
quantum super-$A$-polynomial, provides recursion relations for HOMFLY homology theories for each $r$. Further examples and properties of super-$A$-polynomials have been analyzed in \cite{FGSS,Nawata}.

In this note we summarize the construction of the super-$A$-polynomial and illustrate it in several examples, following (mostly)  \cite{superA,FGS}. We start by recalling the original volume conjectures in section \ref{sec-volume}. In section \ref{sec-new} we generalize these conjectures in a way which leads to the super-$A$-polynomial and the quantum super-$A$-polynomial. Essential ingredients which make these new conjectures work are colored superpolynomials, introduced in section \ref{sec-super}; we stress that developments of tools and techniques which enable to derive an explicit form of such colored superpolynomials is an important, independent result of the work reported here. In section \ref{sec:quantizability} we discuss quantizability properties of super-$A$-polynomials, and in section \ref{sec:phys} we present their interpretation in 3d, $\cN=2$ SUSY gauge theories.


\section{Volume conjectures}  \label{sec-volume}

Originally the ``volume conjecture'' referred to the observation \cite{Kashaev} that the so-called Kashaev invariant of a knot $K$
defined at the $n$-th root of unity $q = e^{2 \pi i / n}$ in the classical limit has a nice asymptotic behavior
determined by the hyperbolic volume $\text{Vol} (M)$ of the knot complement $M = S^3 \setminus K$.
Shortly after, it was realized \cite{MurMur} that the Kashaev invariant is equal to the $n$-colored Jones polynomial of a knot $K$
evaluated at $q = e^{2 \pi i / n}$, so that the volume conjecture could be stated simply as
\be
\lim_{n \to \infty} \frac{2 \pi \log |J_n (K; q = e^{2 \pi i / n})|}{n} \; = \; \text{Vol} (M) \,.
\label{VCbasic}
\ee

The physical interpretation of the volume conjecture was proposed in \cite{Apol}.
Besides explaining the original observation \eqref{VCbasic} it immediately led to a number
of generalizations, in which the right-hand side is replaced by a function of various parameters (see \cite{DGreview} for a review).

\subsection{Generalized volume conjecture}

Once the volume conjecture is put in the context of analytically continued Chern-Simons theory,
it becomes clear that the right-hand side is simply the value of the classical $SL(2,\C)$ Chern-Simons action functional on a knot complement $M$.
Since classical solutions in Chern-Simons theory ({\it i.e.} flat connections on $M$) come in families,
parametrized by the holonomy of the gauge connection on a small loop around the knot,
this physical interpretation immediately leads to a ``family version'' of the volume conjecture \cite{Apol}:
\be
J_n (K; q = e^{\hbar}) \;\overset{{n \to \infty \atop \hbar \to 0}}{\sim}\;
\exp \left( \frac{1}{\hbar} S_0 (u) \,+\, \ldots 
\right)
\label{VCparam}
\ee
parametrized by a complex variable $u$. Here, the limit on the left-hand side is slightly more interesting
than in \eqref{VCbasic} and, in particular, also depends on the value of the parameter $u$:
\be
q = e^{\hbar} \to 1 \,, \qquad
n \to \infty \,, \qquad
q^n = e^u \equiv x  ~~~\text{(fixed)}
\label{VClimit}
\ee
In fact, Chern-Simons theory predicts all of the subleading terms in the $\hbar$-expansion denoted by ellipsis in \eqref{VCparam}.
These terms are the familiar perturbative coefficients of the $SL(2,\C)$ Chern-Simons partition function on $M$.

\subsection{Quantum volume conjecture}

Classical solutions in Chern-Simons theory ({\it i.e.} flat connections on $M$) are labeled
by the holonomy eigenvalue $x = e^u$ or, to be more precise, by a point on the algebraic curve
\be
\cC: \quad \left\{(x,y)\in\mathbb{C}^*\times \mathbb{C}^*\Big|
A(x,y) \; = \; 0 \right\}\,,
\label{Acurve}
\ee
defined by the zero locus of the $A$-polynomial, a certain classical invariant of a knot.
In quantum theory, $A(x,y)$ becomes an operator $\hat A (\hat x, \hat y; q)$
and the classical condition \eqref{Acurve} turns into a statement that the Chern-Simons partition function
is annihilated by $\hat A (\hat x, \hat y; q)$. This statement applies equally well to Chern-Simons theory
with the compact gauge group $SU(2)$ that computes the colored Jones polynomial $J_n (K;q)$ as well as to
its analytic continuation that localizes on $SL(2,\C)$ flat connections.
In the former case, one arrives at the ``quantum version'' of the volume conjecture \cite{Apol}:
\be
\hat A \; J_* (K;q) \; \simeq \; 0 \,,
\label{VCquant}
\ee
which in the mathematical literature was independently proposed around the same time~\cite{Garoufalidis} and is known as the AJ-conjecture.
The action of the operators $\hat x$ and $\hat y$ follows from quantization of Chern-Simons theory, and 
one finds that $\hat x$ acts as a multiplication by $q^n$, whereas $\hat y$ shifts the value of~$n$:
\begin{align}
& \hat x J_n \; = \; q^n J_n \label{xyactionJ} \\
& \hat y J_n \; = \; J_{n+1} \notag
\end{align}
In particular, one can easily verify that these operations obey the commutation relation
\be
\hat y \hat x \; = \; q \hat x \hat y
\label{xycomm}
\ee
that follows from the symplectic structure on the phase space of Chern-Simons theory.
Therefore, upon quantization a classical polynomial relation of the form \eqref{Acurve} becomes
a $q$-difference equation for the colored Jones polynomial or Chern-Simons partition function.
Further details, generalizations, and references can be found in \cite{DGreview}.



\section{New volume conjectures and the super-A-polynomial}   \label{sec-new}

Besides the ``non-commutative'' deformation \eqref{VCquant},
the $A$-polynomial also admits two commutative deformations that
in a similar way encode the ``color behavior'' of two natural generalizations of the colored Jones polynomial:
the $t$-deformation that corresponds to the categorification of colored Jones invariants \cite{FGS}
and $Q$-deformation that corresponds to extending $J_n (K;q)$ to higher rank knot polynomials \cite{superA,AVqdef}.
Our task in what follows is to combine these two deformations into a single unifying structure. In particular, this leads
to a new unifying knot invariant. We call this invariant the {\it super-$A$-polynomial} since it
describes how the $S^{n-1}$-colored superpolynomials $\P_n (a,q,t)$ depend on color, {\it i.e.} on the representation $R = S^{n-1}$,
much in the same way as $A$-polynomial does it for the colored Jones polynomial.

We recall that, in the context of BPS states, the superpolynomial
is defined as a generating function of refined open BPS invariants
on a rigid Calabi-Yau 3-fold $X$ in the presence of a Lagrangian brane supported on $L \subset X$:
\be
\P (a,q,t) \; := \; \Tr_{\cH^{\text{ref}}_{\text{BPS}}} \; a^{\beta} q^{P} t^{F} \,,
\qquad \beta \in H_2 (X,L)
\label{Paqt}
\ee
and, in application to knots, the superpolynomial $\P (K;a,q,t)$ is defined as a Poincar\'e polynomial
of the triply-graded homology theory $\cH (K)$ that categorifies the HOMFLY polynomial $P(K;a,q)$, see \cite{DGR} for details.
According to the conjecture of \cite{GSV}, these two definitions give the same result
when $X$ is the total space of the $\cO (-1) \oplus \cO (-1)$ bundle over $\cp^1$
and $L$ is the Lagrangian submanifold determined by the knot $K \subset \mathbf{S}^3$, {\it cf.} \cite{OoguriV,Taubes,Koshkin}.
Lagrangian branes of multiplicity $r=n-1$ yield the so-called ``$n$-colored'' version
of the superpolynomial which, in the context of knot homologies, was recently introduced in \cite{GS},
\be
\P_n (K; a,q,t) \; := \; \sum_{i,j,k} \, a^i q^j t^k \, \dim \cH^{S^{n-1}}_{i,j,k} (K) \,,
\label{superPdef}
\ee
as a Poincar\'e polynomial of a triply-graded homology theory categorifying
the $S^r$-colored HOMFLY polynomial (see also \cite{FGS,MMSS,ItoyamaMMM}).
For $t=-1$ the above expression reduces to the Euler characteristic and reproduces normalized (i.e. such that $P_n(\unknot;a,q)=1$ for the unknot $\unknot$) colored HOMFLY polynomial $P_n(K;a,q)$:
\be
P_n(K;a,q) = \P_n (K; a,q,-1) = \sum_{i,j,k} \, a^i q^j (-1)^k \, \dim \cH^{S^{n-1}}_{i,j,k} (K) \,.
\label{cHOMFLY}
\ee

\begin{figure}[ht]
\centering
\includegraphics[width=3.0in]{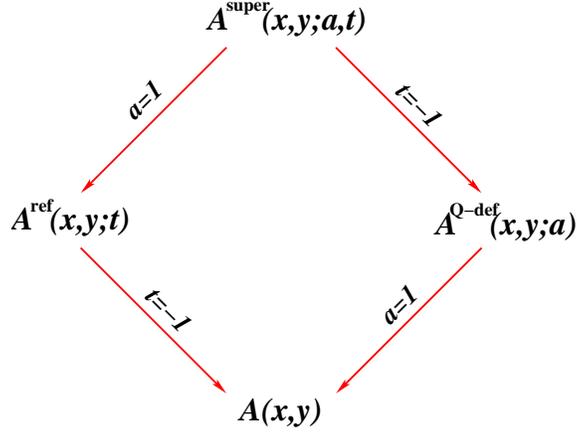}
\caption{Various specializations of the super-$A$-polynomial.}
\label{fig:superAlimits}
\end{figure}

Our main goal is to explain that $S^{n-1}$-colored superpolynomials $\P_n (K;a,q,t)$
depend on color ({\it i.e.} on the representation $R = S^{n-1}$) in a simple and controllable way,
governed by the super-$A$-polynomial $A^{\text{super}} (x,y;a,t)$
and by its quantization $\hat A^{\text{super}} (\hat x, \hat y;a,q,t)$.
Specifically, based on the physics arguments and the study of examples, we propose the following
analog of the generalized volume conjecture \cite{Apol} or its refined version~\cite{FGS}:

\medskip
\noindent
{\bf Conjecture 1:}
{\it In the limit}
\be
q = e^{\hbar} \to 1 \,, \qquad a = \text{fixed} \,, \qquad t = \text{fixed} \,, \qquad x = q^n = \text{fixed}
\label{reflimit}
\ee
{\it the $n$-colored superpolynomials $\P_n (K;a,q,t)$ exhibit the following ``large color'' behavior:}
\be
\P_n (K;a,q,t) \;\overset{{n \to \infty \atop \hbar \to 0}}{\sim}\;
\exp\left( \frac{1}{\hbar} \int \log y \frac{dx}{x} \,+\, \ldots \right)
\label{VCsuper}
\ee
{\it where ellipsis stand for regular terms (as $\hbar \to 0$) and the leading term is given by
the integral on the zero locus of the super-$A$-polynomial, {\it cf.} \eqref{Acurve}:}
\be
A^{\text{super}} (x,y;a,t) \; = \; 0 \,.
\label{supercurve}
\ee

Moreover, just like the ordinary $A$-polynomial has its quantum analog \eqref{VCquant},
the super-$A$-polynomial is a characteristic polynomial of a quantum operator $\hat A^{\text{super}} (\hat x, \hat y;a,q,t)$
that combines commutative $t$- and $a$-deformations with the non-commutative $q$-deformation \eqref{xycomm}.
We call this operator the \emph{quantum super-$A$-polynomial}.

\medskip
\noindent
{\bf Conjecture 2:}
{\it For a given knot $K$, the colored superpolynomial $\P_n (K;a,q,t)$ satisfies a recurrence relation of the form \eqref{VCquant}:}
\be
a_k \, \P_{n+k} (K;a,q,t) + \ldots + a_1 \, \P_{n+1} (K;a,q,t) + a_0 \, \P_n (K;a,q,t) \; = \; 0
\label{QVCsuper}
\ee
{\it where $\hat x$ and $\hat y$ act on $\P_n (K;a,q,t)$ as in \eqref{xyactionJ},
and where the rational functions $a_i \equiv a_i (\hat x, a, q, t)$
are the coefficients of the ``quantum super-$A$-polynomial''}
\be
\hat A^{\text{super}} (\hat x, \hat y; a,q,t) \; = \; \sum_i a_i (\hat x, a, q, t) \, \hat y^i \,,
\label{Asuperform}
\ee
{\it whose characteristic polynomial is} $A^{\text{super}} (x,y;a,t)$.

As in \eqref{VCquant}, sometimes we informally write \eqref{QVCsuper} in the compact form
\be
\hat A^{\text{super}} \; \P_* (K;a,q,t) \; = \; 0 \,,   \label{AsuperP}
\ee
which is a quantum version of the classical curve \eqref{supercurve}.

\begin{table}[h]
\centering
\begin{tabular}{|@{$\Bigm|$}c|c@{$\Bigm|$}c|c@{$\Bigm|$}|}
\hline
\rule{0pt}{5mm}
\textbf{Quantum operator} & \textbf{provides recursion for} & \textbf{classical limit}  \\[3pt]
\hline
\hline
\rule{0pt}{5mm}
~~$\hat A^{\text{super}} (\hat x, \hat y; a,q,t)$~~ & ~~colored superpolynomial~~~ & ~~~$A^{\text{super}} (x,y;a,t)$~~  \\[3pt]
\hline
\rule{0pt}{5mm}
$\hat A^{\text{ref}} (\hat x, \hat y; q,t)$ & colored $sl(2)$ homology & $A^{\text{ref}} (x,y;t)$  \\[3pt]
\hline
\rule{0pt}{5mm}
$\hat A^{\text{Q-def}} (\hat x, \hat y; a,q)$ & colored HOMFLY & $A^{\text{Q-def}} (x,y;a)$  \\[3pt]
\hline
\rule{0pt}{5mm}
$\hat A (\hat x, \hat y; q)$ & colored Jones & $A (x,y)$  \\[3pt]
\hline
\end{tabular}
\caption{Quantum super-$A$-polynomial and its specializations lead to recursion relations for various $S^n$-colored knot invariants.\label{AAAtable}}
\end{table}

The superpolynomial unifies many polynomial and homological invariants of knots
that can be obtained from it via various specializations, applying differentials, {\it etc.}
For example, for $\cH$-thin knots the specialization to $a=q^2$ yields
the Poincar\'e polynomial of the colored $sl(2)$ knot homology.
Therefore, if $K$ is a thin knot ({\it e.g.} if $K$ is a two-bridge knot),
in the limit $a=q^2$ we expect \eqref{VCsuper} and \eqref{QVCsuper} to reproduce
the corresponding versions of the refined volume conjectures proposed in~\cite{FGS}.
In particular,
\be
\hat A^{\text{super}} (\hat x, \hat y; a=q^2,q,t) \; = \; \hat A^{\text{ref}} (\hat x, \hat y; q,t) \,,
\label{AsuperAref}
\ee
and, via further specialization to the classical limit $q=1$,
\be
A^{\text{super}} (x, y; a=1,t) \; = \; A^{\text{ref}} (x, y; t) \,.
\ee

Similarly, the specialization of the superpolynomial $\P_n (K;a,q,t)$
to $t=-1$ yields the HOMFLY polynomial or, in the problem at hand,
the {\it colored} HOMFLY polynomial~\cite{GS}. Therefore, at $t=-1$ the recursion relation \eqref{QVCsuper}
should reduce to the recursion relation for the $S^{n-1}$-colored HOMFLY polynomial,
whose characteristic variety --- called the $Q$-deformed $A$-polynomial in \cite{AVqdef} --- must be
contained in $A^{\text{super}} (x,y;a,t=-1)$ as a factor.
To avoid clutter, we include possible extra factors inherited from $A^{\text{super}} (x,y;a,t)$
in the definition of the $Q$-deformed $A$-polynomial, so that
\be
A^{\text{super}} (x, y; a,t = -1) \; = \; A^{\text{Q-def}} (x, y; a) \,.
\label{AsuperAaug}
\ee
Moreover, the authors of \cite{AVqdef} proposed an important conjecture that offers a new way of
looking at this polynomial (that, in our Figure \ref{fig:superAlimits}, occupies the right corner)
and identifies it with the augmentation polynomial of knot contact homology~\cite{NgFramed}.
In what follows we use the names ``$Q$-deformed $A$-polynomial''
and ``augmentation polynomial'' interchangeably.
In fact, one justification for this comes from the fact (see \cite[Proposition 5.9]{NgFramed} for a proof)
that the classical augmentation polynomial, when specialized further to $a=1$,
reduces to the ordinary $A$-polynomial,
possibly with some extra factors, which altogether we denote simply by $A(x,y)$:
\be
A^{\text{super}} (x,y; a=1,t = -1) \; = \; A^{\text{Q-def}} (x,y; a=1) \; = \; A (x,y) \,,
\label{AugA}
\ee
as it should in order to fit perfectly in the diagram in Figure \ref{fig:superAlimits}.

Therefore, our super-$A$-polynomial $A^{\text{super}} (x,y;a,t)$ can be viewed,
on one hand, as a ``refinement'' of the augmentation polynomial $A^{\text{Q-def}} (x,y;a)$
and, on the other hand, as a ``$Q$-deformation'' of the refined $A$-polynomial $A^{\text{ref}} (x,y;t)$,
see Figure \ref{fig:superAlimits}.
For important examples of super-$A$-polynomials (to be discussed in more detail in what follows) see table~\ref{table:superA}.
Other interesting specializations of super-$A$-polynomials, \emph{e.g.} involving setting $x=1$ or $q=1$, are discussed in \cite{superA,FGSS}.

\begin{table}[h]
\centering
\begin{tabular}{|@{$\Bigm|$}c|@{$\Bigm|$}l|}
\hline
\textbf{Knot}  & $A^{\text{super}} (x,y;a,t)$  \\
\hline
\hline
Unknot, $\unknot$ & {\scriptsize $ ( - a^{-1} t^{-3})^{1/2} (1 + a t^3 x)  - (1 - x) y $}  \\
\hline
 & {\scriptsize $ a^2 t^5 (x-1)^2 x^2 + a t^2 x^2 (1 + a t^3 x)^2 y^3 +$}  \\
Figure-eight, ${\bf 4_1}$ & {\scriptsize $+ a t (x-1) (1 + t(1-t) x  + 2 a t^3(t+1) x^2
    -2 a t^4(t+1) x^3  + a^2 t^6(1-t) x^4  - a^2 t^8 x^5) y$} \\
    & {\scriptsize $- (1 + a t^3 x) (1 + a t(1-t) x +
    2 a t^2(t+1) x^2  + 2 a^2 t^4(t+1) x^3  + a^2 t^5(t-1) x^4  + a^3 t^7 x^5) y^2 $} \\
\hline
Trefoil, ${\bf 3_1}$ & {\scriptsize $a^2 t^4 (x-1) x^3 -a\big( 1 - t^2 x + 2 t^2 (1 + a t) x^2 + a t^5 x^3 + a^2 t^6 x^4 \big) y + (1 + a t^3 x) y^2$}  \\
\hline
$(2,2p$+$1)$ torus knot & {\scriptsize eliminate $z_0$ in   $ \left\{\begin{array}{l} 1 = \frac{(z_0-x)(t^2z_0-1)(1+ at^3 xz_0)}{t^{2+2p}z_0^{1+2p}(z_0-1)(atx+z_0)(t^2 x z_0-1)}   \\
y =  \frac{a^p t^{2 + 2 p} (x-1) x^{1 + 2 p} (atx + z_0) (1 + a t^3 x z_0)}{(1 + a t^3 x) (x - z_0) (t^2 x z_0-1)}   \end{array} \right.\  $ see table \protect\ref{table_super}  }  \\
\hline
\end{tabular}
\caption{Super-$A$-polynomials for simple knots.\label{table:superA}}
\end{table}


\section{Essential ingredients: colored superpolynomials}    \label{sec-super}

The (quantum) super-$A$-polynomial arising from the conjectures presented in section \ref{sec-new} is intimately related to the colored superpolynomials introduced in (\ref{superPdef}). Indeed, the knowledge of $S^{n-1}$-colored superpolynomials $\P_n (K; a,q,t)$ for general $n$ allows to determine the (quantum) super-$A$-polynomial, and this is how super-$A$-polynomials will be derived in all examples in what follows. Nonetheless, determining $S^{n-1}$-colored superpolynomials is itself a hard task, and these objects are not even defined mathematically in a rigorous and computable way. However it turns out that physics offers two possible ways to obtain (or, at least, to conjecture) the form of colored superpolynomials: either using the so-called refined Chern-Simons theory, or taking advantage of the structure of differentials in homological theories. For some knots both of these methods can be used, and then they lead to remarkable identities, which confirm validity of the physics approach. We need to get acquainted with these methods before we present the construction and examples of super-$A$-polynomials.

\subsection{Refined Chern-Simons theory}

As is well known \cite{Witten_Jones}, knot invariants (more precisely -- quantum group invariants associated to a given knot) are simply related to expectation values of Wilson loop operators $W_R(K)[A]:={\rm Tr}_R \text{P} \exp\left[\oint_K A\right]$, supported on a knot $K$ and decorated by a representation $R$, in Chern-Simons theory
\be
Z_{G}^{\text{CS}}(M,K_R;q)
\; = \; \int [dA] \, W_R(K)[A] \, e^{ik S_{\text{CS}}[A;M]},
\label{ZCSdef}
\ee
where the Chern-Simons action on a 3-manifold $M$ reads
\begin{eqnarray}
S_{\text{CS}}[A;M] \; = \; \frac{1}{4\pi}\int_M {\rm Tr}_{\text{adj}}
\left( A\wedge dA+\frac{2}{3}A\wedge A\wedge A\right) \,.
\label{SCSdef}
\end{eqnarray}
The quantum group invariants are reproduced as the above expectation values normalized by that of the unknot (which we often denote as~$\unknot$), and remarkably such expressions are simply polynomials in $q=e^{\frac{2\pi i}{k+h}}$ with integer coefficients (at least when $M=\S^3$)
\begin{eqnarray}
J^{\frak{g}, R} (K;q) \; =\; \frac{Z_{G}^{\rm CS}(\S^3,K_R;q)}{Z_{G}^{\rm CS}(\S^3,\unknot_R;q)} \,.
\label{JasZZratio}
\end{eqnarray}
In particular, for $\frak{g} = sl(N)$ a dependence on $N$ is very simple and $J^{sl(N), R} (K;q)$ turn out to be polynomials in $q$ and $a=q^N$ which reproduce (normalized) colored HOMFLY polynomials $J^{sl(N), R} (K;q) = P^R (K; a = q^N, q)$. For $R=S^{n-1}$ they are denoted $P_n (K; a, q)$, and they already appeared above in (\ref{cHOMFLY}). 


Therefore our task is to introduce, from the perspective of the Chern-Simons theory, a dependence on the Poincar\'e variable $t$ into the colored HOMFLY polynomial, so that Chern-Simons amplitudes given in (\ref{JasZZratio}) would be extended to $t$-dependent quantities, which should be identified with (\ref{superPdef}). Such a generalization has been introduced in \cite{AS} and is often referred to as refined Chern-Simons theory. More precisely, the fundamental formulation (in terms of a $t$-dependent action) of refined Chern-Simons theory is still not known. Nonetheless, the authors of \cite{AS} argued how a dependence on $t$ should be introduced in a consistent manner in various quantities arising in the quantization of the original Chern-Simons theory. In particular they proposed a refined version of modular matrices $S$ and $T$ which satisfy the Verlinde formula. Similarly other objects in Chern-Simons theory become functions of $q$ and $t$, or equivalently (as often arising in various calculations) $q_1$ and $q_2$ (such that $q = \frac{1}{q_2},  t = - \sqrt{\frac{q_2}{q_1}} $). For example, Schur polynomials that arise in original Chern-Simons theory (in particular as expectation values of the unknot) are replaced by Macdonald polynomials, and so on:
\begin{eqnarray}
\underline{\text{~~CS gauge theory~~}} && \underline{\text{~~refined invariants~~}}
\nonumber \\
Z_{SU(N)}^{\rm CS}(\S^3,K_R;q)
&\qquad \leadsto \qquad&
Z^{\text{ref}}_{SU(N)}(\S^3,K_R;q_1,q_2)
\nonumber \\
\dim_q R = s_R(q^{\varrho})
&\qquad \leadsto \qquad&
M_R(q_2^{\varrho};q_1,q_2)
\nonumber \\
q^{C_2(R)}
&\qquad \leadsto \qquad&
q_1^{\frac{1}{2}||R||^2}q_2^{-\frac{1}{2}||R^{t}||^2}q_2^{\frac{N}{2}|R|}q_1^{-\frac{1}{2N}|R|^2}
\label{Casimir_ref} \\
& \vdots & \nonumber
\end{eqnarray}
Refined Chern-Simons theory is still quite mysterious; in particular explicit computations are possible only for some particular knots, and they involve various subtleties, related \emph{e.g.} to the appearance of the so-called $\gamma$-factors (which are irrelevant for $t=-1$). Nonetheless, we are able to predict an explicit form of superpolynomials using refined Chern-Simons theory in various examples, such as the (unnormalized) unknot (see (\ref{Punknot})), or $(2,2p+1)$ torus knots (for arbitrary $p$, see (\ref{Paqt-torus})). 
Various aspects of refined Chern-Simons theory are discussed in detail in \cite{AS,FGS,DMMSS}.


\subsection{Differentials in knot homologies}   

Even though combinatorial definition of colored knot homologies (\ref{superPdef}) is, in general, not known, it turns out that various physics arguments predict how the action of various differentials in knot homologies should look like. Such differentials endow knot homologies with a very rich structure,
which turns out to be very elegant and often so constraining that one can even compute
colored superpolynomials $\P_n (K; a,q,t)$ based on this structure alone, with a minimal input.
In particular, this is how nice formulas like (\ref{Paqt41}), \eqref{Paqt31}, or (\ref{fort2k}) can be produced.
Referring the reader to \cite{GS,FGSS} for further details, here we merely state
a simple rule of thumb: the factors of the form $(1 + a^i q^j t^k)$ that we often
see {\it e.g.} in
\eqref{Punrec},
\eqref{Paqt41},
\eqref{Paqt31}, and \eqref{Paqt-torus}
come from differentials of $(a,q,t)$-degree $(i,j,k)$, {\it cf.}~\cite[eq. (3.54)]{FGS}:
\be
\begin{array}{c@{\;}|@{\;}c@{\;}|@{\;}c@{\;}c}
\text{differentials} & \text{factors} & (a, q, t)~\text{grading} \\\hline
d_{N>0} & \quad 1 + a q^{-N} t \quad & (-1,N,-1) \\[.1cm]
d_{N<0} & 1 + a q^{-N} t^3 & (-1,N,-3) \\[.1cm]
d_{\text{colored}} & 1 + q & (0,1,0) \\[.1cm]
& 1 + at & (-1,0,-1) \\
& \vdots &
\end{array}
\label{gradingtabl}
\ee
For example, notice that all terms with $k>0$ in the expression \eqref{Paqt41}
for the colored superpolynomial of the figure-eight knot
manifestly contain a factor $(1 + a q^{n-1} t^3)$.
Hence, the $S^{n-1}$-colored superpolynomial of the figure-eight knot
has the following structure 
\be
\P_n ({\bf 4_1};a,q,t) \; = \; 1 + (1 + a q^{n-1} t^3) Q_n (a,q,t) \,,
\ee
which means that, when evaluated at $a = - q^{1-n} t^{-3}$,
the sum \eqref{Paqt41} collapses to a single $k=0$ term,
$\P_n ({\bf 4_1}; a = - q^{1-n} t^{-3},q,t) = 1$.
A proper interpretation of this fact is that a specialization to $N = 1-n$
of the triply-graded $S^{n-1}$-colored HOMFLY homology, carried out by the action of the differential $d_{1-n}$, is trivial.
In other words, the differential $d_N$ with $N=1-n$ is canceling in a theory with $R = S^{n-1}$.
A systematic implementation of such observations fully determines the form of colored superpolynomials  -- at least for some knots -- as we will see in the examples in the next section.


\section{Examples}    \label{sec-case}

In this section we illustrate ideas presented above in explicit examples of various knots. We start with the simplest example of the unknot, and then discuss a non-trivial example of a hyperbolic knot, \emph{i.e.} the figure-eight knot, and the entire family of $(2,2p+1)$ torus knots, with a special emphasis on the trefoil. In each case we start our considerations by providing explicit and general formulas for $S^{n-1}$-colored superpolynomials $\P_n (a,q,t)$, illustrating the power of two approaches described in section \ref{sec-super}. Taking advantage of these representations, subsequently we derive classical and quantum super-$A$-polynomials for these knots and discuss their properties. For other examples of superpolynomials and super-$A$-polynomials see \cite{FGSS,Nawata}.

\subsection{Unknot}
\label{sec:unknot}

Let us start with the simplest example of the unknot. Despite its simplicity, this is still an interesting and important example; as we will see, some objects associated to the unknot, which are trivial in the non-refined and non-super case, become rather non-trivial when $t$- or $a$-dependence is turned on.

We recall than in the unknot case we must consider unreduced (or ``unnormalized'') knot polynomials --
in particular, unreduced colored superpolynomial $\bar{\P}_{n}(a,q,t)$ -- since, by definition,
reduced polynomials are normalized by the value of the unknot, so that $\P_{n}(\unknot; a,q,t)=1$.
From the viewpoint of the (refined) Chern-Simons theory the unreduced
colored superpolynomial is defined as the ratio of partition functions
on $\S^3$ in the presence and absence of a knot. In case of the unknot
this ratio is given by the Macdonald polynomial, and after the change of
variables $q = \frac{1}{q_2},  t = - \sqrt{\frac{q_2}{q_1}} $, 
we find that the $S^{n-1}$-colored superpolynomial reads
\begin{eqnarray}
\bar{\P}_n (\unknot; a,q,t)
&=& \frac{Z_{SU(N)}^{\text{ref}}(\S^3,\unknot_{\Lambda^{n-1}};q_1,q_2)}{Z_{SU(N)}^{\text{ref}}(\S^3;q_1,q_2)}
=M_{\Lambda^{n-1}}(q_2^{\varrho};q_1,q_2) \nonumber \\
&=&
(-1)^{\frac{n-1}{2}}a^{-\frac{n-1}{2}}q^{\frac{n-1}{2}}t^{-\frac{3(n-1)}{2}} \frac{(-at^3;q)_{n-1}}{(q;q)_{n-1}} \,.
\label{Punknot}  
\end{eqnarray}
We also recall that the $q$-Pochhammer symbol $(x,q)_k$, defined by a product formula, has the following asymptotics
\be
(x,q)_k \equiv \prod_{i=0}^{k-1} (1 - x q^i) \sim e^{\frac{1}{\hbar}\left({\rm Li}_2(x)-{\rm Li}_2(x q^k)\right)}  . \label{qPoch}
\ee

\begin{figure}[ht]
\begin{center}
\includegraphics[width=0.6\textwidth]{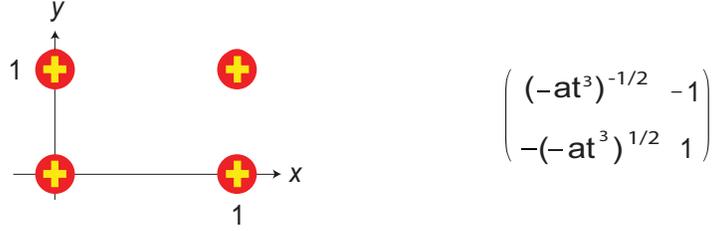}
\end{center}
\caption{Newton polygon for the super-$A$-polynomial of the unknot (left).
Red circles denote monomials of the super-$A$-polynomial, and smaller yellow crosses denote monomials of its $a=-t=1$ specialization.
In this example both Newton polygons look the same, so that positions of all circles and crosses overlap.
The coefficients of the super-$A$-polynomial are also shown in the matrix on the right.
The role of rows and columns is exchanged in these two presentations:
a monomial $a_{i,j}x^i y^j$ corresponds to a circle (resp. cross) at position $(i,j)$ in the Newton polygon,
while in the matrix on the right it is shown as the entry $a_{i,j}$ in the $(i+1)^{\text{th}}$ row and in the $(j+1)^{\text{th}}$ column.
These conventions are the same as in~\cite{FGS,superA,FGSS}.}
\label{fig-unknotNewtonMatrix}
\end{figure}

Once the general expression for the colored superpolynomial is determined, we can find a recursion relation it satisfies.
In particular, as the homological unknot invariant (\ref{Punknot}) has a product form,
we can immediately write down the recursion relation it satisfies:
\be
\bar{\P}_{n+1} (\unknot; a,q,t) \; = \; (-a^{-1} t^{-3}q)^{1/2}  \frac{1 + a t^3 q^{n-1}}{1 - q^{n}}  \, \bar{\P}_n (\unknot; q,t) \,.
\label{Punrec}
\ee
This means  that the quantum super-$A$-polynomial for the unknot reads
\be
\hat A^{\text{super}}(\hat{x},\hat{y};a,q,t) \; = \;
( - a^{-1} t^{-3} q )^{1/2} (1 + a t^3 q^{-1} \hat{x})  - (1 - \hat x) \hat y \,.   \label{qAt-unknot}
\ee
In the classical limit $q \to 1$ this operator
reduces to the  classical super-$A$-polynomial 
\be
A^{\text{super}} (x,y;a,t) \; = \;
( - a^{-1} t^{-3})^{1/2} (1 + a t^3 x)  - (1 - x) y \,.   \label{At-unknot}
\ee
The Newton polygon as well as the coefficients of monomials of this polynomial are shown in figure~\ref{fig-unknotNewtonMatrix}.
In the unrefined limit $t=-1$ the relation (\ref{qAt-unknot}) takes the form
\be
\hat A^{\text{Q-def}}(\hat x, \hat y ;a,q) \; = \;
(a^{-1} q)^{1/2}(1 - a q^{-1} \hat x)  - (1 - \hat x) \hat y \,,
\ee
and specializing further to $q=1$ we get the augmentation polynomial
\be
A^{\textrm{Q-def}}(x,y;a) \; = \; a^{-1/2}(1 - a x)-(1-x)y \,.   \label{Aunknot1x1y}
\ee
Interestingly, this polynomial does not factorize, and only in the limit of ordinary $A$-polynomial $a\to 1$ do we get a factorized form with $y-1$ factor representing the abelian connection
\be
A(x,y) \; = \; (1 -  x)(1-y) \,.  \label{Axy-factor}
\ee

It is instructive to show that the super-$A$-polynomial (\ref{At-unknot}) can be also derived from the asymptotic analysis of~(\ref{Punknot}).
Indeed, using the asymptotics (\ref{qPoch}), in the limit (\ref{reflimit}) we can approximate (\ref{Punknot}) as
$$
P_n(\unknot; a,q,t) = \exp \frac{1}{\hbar}\Big(\log x\, \log(-a^{-1} t^{-3})^{1/2} + \textrm{Li}_2(x) - \textrm{Li}_2(-a t^3 x)  + \textrm{Li}_2(-a t^3) - \frac{\pi^2}{6} + \mathcal{O}(\hbar) \Big),
$$
from which identify the potential $\widetilde{\mathcal{W}} = \int \log y \frac{dx}{x}$ in \eqref{VCsuper} as
\be
\widetilde{\mathcal{W}} \; = \;  \log x\, \log(-a^{-1} t^{-3})^{1/2} + \textrm{Li}_2(x) - \textrm{Li}_2(-a t^3 x)  + \textrm{Li}_2(-a t^3) - \frac{\pi^2}{6} \,. \label{S0unknot}
\ee
Differentiating it with respect to $x$, we now obtain
\be
y = e^{x\partial_x \widetilde{\mathcal{W}} } = (-a^{-1} t^{-3})^{1/2} \frac{1 + a t^3 x}{1 - x} \,,
\ee
which reproduces the defining equation of the super-$A$-polynomial given in (\ref{At-unknot}).
We also note that for $a=-t=1$ the potential $\widetilde{\mathcal{W}}$ vanishes,
which is related to the factorization occurring in (\ref{Axy-factor}) and can be attributed
to the fact that the only $SL(2,\C)$ flat connections on a solid torus (= complement of the unknot) are abelian flat connections.
When $a\neq 1$ or $t\neq -1$, the potential $\widetilde{\mathcal{W}}$ is nonzero and presumably can be interpreted
as a contribution of ``deformed'' abelian flat connections. 


\subsection{Figure-eight knot}
\label{sec:fig8}

In this section we consider the figure-eight knot, also denoted ${\bf 4_1}$.
This is a hyperbolic knot, and we stress that it provides a highly non-trivial example,
for which many simplifications common in the realm of torus knots (to be discussed in the following sections) do not occur.

The colored superpolynomial (\ref{Paqt}) for figure-eight knot can be found using the highly constraining structure of differentials. This strategy has been employed in \cite{FGSS} and the resulting superpolynomial reads
\be
\P_n ({\bf 4_1};a,q,t) = \sum_{k=0}^{\infty} (-1)^k a^{-k} t^{-2k} q^{-k(k-3)/2} \frac{(-a t q^{-1},q)_k}{(q,q)_k}  (q^{1-n},q)_k (-a t^3 q^{n-1}, q)_k \,.
\label{Paqt41}
\ee
An independent, though not entirely unrelated derivation that gives the same result has been proposed in \cite{ItoyamaMMM}. Explicit values of $\P_n ({\bf 4_1};a,q,t)$ for low values of $n$ are given in table~\ref{tab-P41}; note that they are all polynomials with positive coefficients, as necessarily expected from (\ref{superPdef}).
We stress that (\ref{Paqt41}) is in itself a very strong result, which illustrates the power of physics methods; to appreciate this fact and to confirm the validity of the above result we note that:
\begin{itemize}
\item for $a=q^2$ and $t=-1$, the formula (\ref{Paqt41}) reduces to the familiar expression for the colored Jones polynomial studied {\it e.g.} in \cite{Habiro,Garoufalidis}:
$$
J_n ({\bf 4_1};q)=\P_n ({\bf 4_1};q^2,q,-1)= 
\sum_{k=0}^{n-1}  q^{n k} (q^{-n-1},q^{-1})_k (q^{-n+1},q)_k
$$
\item for $t=-1$, (\ref{Paqt41}) agrees with the colored HOMFLY polynomial given in the unpublished work \cite{kawagoe}, which was also used in the analysis of \cite{AVqdef} (for precise relation see \cite{superA});
\item for $n=2$ the superpolynomial (\ref{Paqt41}) agrees with the known result given {\it e.g.} in~\cite{DGR} (to match conventions we need to replace $a$ and $q$ in \cite{DGR} respectively by $a^{1/2}$ and $q^{1/2}$);
\item for $n=3$ and $n=4$ the expression (\ref{Paqt41}) reproduces results given in~\cite{GS};
\item for $a=-q^j t^k$ the expression \eqref{Paqt41} correctly reproduces specializations predicted from the colored / canceling differentials with $(a,q,t)$-grading $(-1,j,k)$, see~\cite{GS}.
\end{itemize}

\begin{table}[h]
\centering
\begin{tabular}{|@{$\Bigm|$}c|@{$\Bigm|$}l|}
\hline
$\, n$ & $\qquad \P_n ({\bf 4_1};a,q,t)$  \\
\hline
\hline
$\, 1$ & $1$   \\
\hline
$\, 2$ & $a^{-1}t^{-2} + t^{-1}q^{-1} + 1 + q t + a t^2$
\\
\hline
$\, 3$ & $ a^{-2} q^{-2} t^{-4} + (a^{-1} q^{-3} + a^{-1} q^{-2}) t^{-3} + (q^{-3} + a^{-1} q^{-1} + a^{-1}) t^{-2} +$ \\
& $+ (q^{-2} + q^{-1} + a^{-1} + a^{-1} q) t^{-1} + (q^{-1} + 3 + q) + (q^2 + q + a + a q^{-1}) t +$\\
& $+ (q^3 + a q + a) t^2 + (a q^3 + a q^2) t^3 + a^2 q^2 t^4 $\\
\hline
$\, 4$ & $1 + (1 + a^{-1} q t^{-1}) (1 + a^{-1} t^{-1}) (1 + a^{-1} q^{-1} t^{-1}) \times $\\
 & $\qquad \times   (1 + a^{-1} q^{-3} t^{-3}) (1 + a^{-1} q^{-4} t^{-3}) (1 + a^{-1} q^{-5} t^{-3}) a^3 q^6 t^6 +$\\
 & $+ (1 + q + q^2) (1 + a^{-1} q t^{-1}) (1 + a^{-1} q^{-3} t^{-3}) a t^2 +$ \\
 & $ + (1 + q + q^2) (1 + a^{-1} q t^{-1}) (1 + a^{-1} t^{-1}) (1 + a^{-1} q^{-3} t^{-3}) (1 + a^{-1} q^{-4} t^{-3}) a^2 q^2 t^4$\\
\hline
\end{tabular}
\caption{The colored superpolynomial of the ${\bf 4_1}$ knot for $n=1,2,3,4$.
\label{tab-P41} }
\end{table}

\begin{figure}[ht]
\begin{center}
\includegraphics[width=0.8\textwidth]{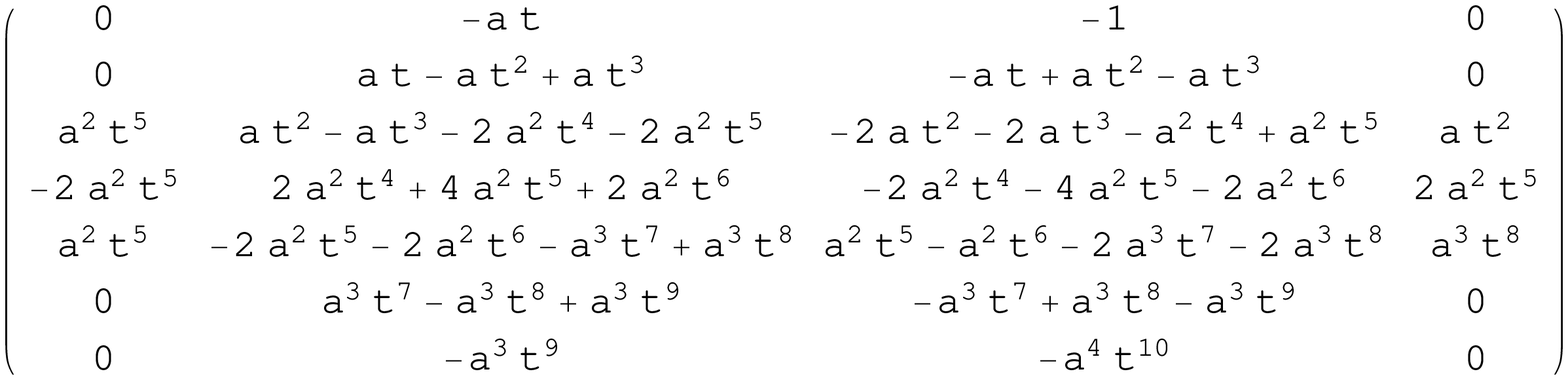}
\caption{Matrix form of the super-$A$-polynomial for the figure-eight knot. The conventions are the same as in the unknot example in figure \protect\ref{fig-unknotNewtonMatrix}.}
\label{fig:matrix41}
\end{center}
\end{figure}

\begin{figure}[ht]
\begin{center}
\includegraphics[width=0.5\textwidth]{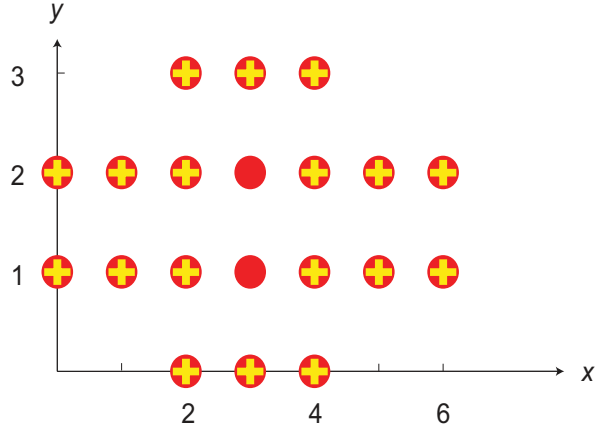}
\caption{Newton polygon of the super-$A$-polynomial for the figure-eight knot and its $a=-t=1$ limit. The conventions are the same as in figure \protect\ref{fig-unknotNewtonMatrix}.}
\label{fig:Newton41}
\end{center}
\end{figure}

Recursion relations satisfied by (\ref{Paqt41}) can be found using the Mathematica package \textrm{qZeil.m} developed by \cite{qZeil}.
In the notation of (\ref{Asuperform}) these recursions take form
\be
\hat A^{\text{super}} (\hat x, \hat y; a,q,t) = a_0 + a_1 \hat{y} + a_2 \hat{y}^2 + a_3 \hat{y}^3,
\ee
where
\bea
a_0 & = &  \frac{a t^3 (1 - \hat{x}) (1 - q \hat{x}) (1 + a t^3 q^2 \hat{x}^2) (1 +
   a t^3 q^3 \hat{x}^2 )}{q^3 (1 + a t^3 \hat{x}) (1 + a  t^3 \hat{x}^2) (1 +
   a t^3 q \hat{x} ) (1 + a t^3 q^{-1} \hat{x}^2)}  \nonumber \\
a_1 & = & - \frac{(1 - q \hat{x}) (1 + a t^3 q^3\hat{x}^2 )}
 {t q^3 \hat{x}^2  (1 + a t^3 \hat{x}) (1 + a t^3 q \hat{x}) (1 + a t^3 q^{-1} \hat{x}^2)}   \nonumber \\
& & \quad \times \Big( 1 - t (t-1) q \hat{x} + a t^3 q^{-1}(1 + q^3 + q t + q^2 t) \hat{x}^2    \nonumber \\
 & & \qquad - a t^4 (q + q^2 + t + q^3 t)\hat{x}^3  - a^2 (t-1) t^6 q \hat{x}^4 - a^2 t^8 q^2 \hat{x}^5 \Big)   \nonumber \\
a_2 & = &  - \frac{(1 + a t^3 q^2 \hat{x}^2)}{a t^2 q^2 \hat{x}^2 (1 + a t^3 \hat{x}^2) (1 + a t^3 q \hat{x})}   \nonumber \\
& & \quad \times \Big( 1 - a t(t-1) \hat{x} + a  t^2 (q + q^2 + t + q^3 t) \hat{x}^2     \nonumber \\
 & & \qquad  + a^2  t^4 (1 + q^3 + q t + q^2 t) \hat{x}^3 + a^2 (t-1) t^5 q^3 \hat{x}^4 + a^3 t^7 q^3 \hat{x}^5  \Big)   \nonumber \\
a_3 & = &  1   \nonumber
\eea
Taking the classical limit $q\to 1$ (and clearing the denominators), we find the following classical super-$A$-polynomial
\bea
& & A^{\text{super}} (x, y; a,t) \, = \,
a^2 t^5 (x-1)^2 x^2 + a t^2 x^2 (1 + a t^3 x)^2 y^3 + \label{Asuper41} \\
& & \quad + a t (x-1) (1 + t(1-t) x  + 2 a t^3(t+1) x^2
    -2 a t^4(t+1) x^3  + a^2 t^6(1-t) x^4  - a^2 t^8 x^5) y  \nonumber \\
 & & \quad   - (1 + a t^3 x) (1 + a t(1-t) x +
    2 a t^2(t+1) x^2  + 2 a^2 t^4(t+1) x^3  + a^2 t^5(t-1) x^4  + a^3 t^7 x^5) y^2 . \nonumber
\eea
The coefficients of the monomials in this polynomial are assembled into a matrix form presented in figure~\ref{fig:matrix41},
and the corresponding Newton polygon is given in figure~\ref{fig:Newton41}.

According to the Conjecture~1, we should be able to reproduce the same polynomial from the asymptotic behavior of the colored superpolynomial~(\ref{Paqt41}).
This is indeed the case. To show this, we introduce the variable $z=e^{\hbar k}$. Then, in the limit (\ref{reflimit}) with $z=\text{const}$
the sum over $k$ in (\ref{Paqt41}) can be approximated by the integral
\be
\P_n ({\bf 4_1};a,q,t) \; \sim \; \int dz\; e^{\frac{1}{\hbar}\left(\widetilde{\mathcal{W}}({\bf 4_1};z,x)+{\cal O}(\hbar)\right)} \,.   \label{Pn41-integral}
\ee
The potential $\widetilde{\mathcal{W}}({\bf 4_1};z,x)$ can be determined from the asymptotics (\ref{qPoch}):
\bea
& & \widetilde{\mathcal{W}}({\bf 4_1};z,x)  =  \pi i \log z - \frac{\pi^2}{6}- (\log a + 2\log t) \log z - \frac{1}{2} (\log z)^2   \label{V41}\\
& & \quad   + \Li_2( x^{-1}) - \Li_2(x^{-1}z) + \Li_2(-a t) - \Li_2(-a t z) + \Li_2(-a x t^3)
-  \Li_2(-a x t^3 z)  - \Li_2(z) \,.  \nonumber
\eea
At the saddle point
\be
\frac{\partial \widetilde{\mathcal{W}}({\bf 4_1};z,x)}{\partial z}\Bigg|_{z=z_0}=0
\label{saddle_point41}
\ee
it determines the leading asymptotic behavior (\ref{VCsuper}), which at the same time is also computed by the integral along the curve~(\ref{supercurve}),
implying the key identity
\be
y = \exp\left(x\frac{\partial \widetilde{\mathcal{W}}({\bf 4_1};z_0,x)}{\partial x}\right) \,.    \label{yV}
\ee
Plugging the expression (\ref{V41}) to the above two equations we obtain the following system
\be
\left\{\begin{array}{l} 1 = \frac{(x - z_0) (1 + a t z_0) (1 + a t^3 x z_0)}{a t^2 x z_0 (z_0-1)}   \\
y =  \frac{(x-1) (1 + a t^3 x z_0)}{(1 + a t^3 x) (x - z_0)}   \end{array} \right.
\ee
Eliminating $z_0$ from these two equations we indeed reproduce the super-$A$-polynomial~(\ref{Asuper41}).
Overall, the above statements verify the validity of the Conjecture~1 and Conjecture~2 for the figure-eight knot.

Note that for $t=-1$ and $a=1$, the expression (\ref{Asuper41}) reduces to
\be
A(x, y) = (x-1)^2(y-1)\Big(x^2(y^2 + 1) - (1 - x - 2 x^2 - x^3 + x^4) y \Big),  \label{Afigure8}
\ee
which, apart from the $(x-1)^2$ factor, reproduces the $A$-polynomial of the ${\bf 4_1}$ knot,
including the $(y-1)$ factor representing the contribution of abelian flat connections.
We stress that both the factorization and the explicit form of this abelian branch is seen only in the limit $a=-t=1$
and is completely ``mixed'' with the other branches otherwise.
In general the super-$A$-polynomial (\ref{Asuper41}) does not factorize,
as is also the case for the unknot and torus knots that will be discussed next.

More generally, after a simple change of variables
\be
Q = a,\qquad \beta = x, \qquad \alpha = y\frac{1 - \beta Q}{Q(1-\beta)},
\label{super2av_fig8}
\ee
and for $t=-1$ we find that (\ref{Asuper41}) becomes
\bea
A^{\textrm{Q-def}}(\alpha,\beta,Q) & = & \frac{Q^2(1-\beta)^2}{\beta Q - 1}\, \Big( (\beta^2 - Q \beta^3) + (2 \beta - 2 Q^2 \beta^4 + Q^2 \beta^5 - 1) \alpha + \nonumber \\
& & + (1 - 2 Q \beta + 2 Q^2 \beta^4 - Q^3 \beta^5) \alpha^2 + Q^2 (\beta-1) \beta^2 \alpha^3 \Big) \,. \nonumber
\eea
Up to the first fraction, the expression in the big bracket reproduces the $Q$-deformed $A$-polynomial given in \cite{AVqdef}. A related change of variables (for details see \cite{superA}) reveals the relation to the augmentation polynomial of \cite{NgFramed}.


\subsection{Trefoil knot}

In this section, we derive the classical and quantum super-$A$-polynomial for the trefoil knot (\emph{i.e.} $(2,3)$ torus knot, also denoted $T^{(2,3)}$ or ${\bf 3_1}$) and verify the validity of the Conjecture 1 and 2 for this knot.
The analysis follows the same lines as in previous sections, and its starting point is the expression for the colored superpolynomial.
We can provide such an expression from two sources. First, the colored superpolynomial
for general $(2,2p+1)$ torus knot was derived in \cite{FGS} from the perspective of the refined Chern-Simons theory.
This superpolynomial is given in~(\ref{Paqt-torus}), as we will need it for the analysis of general torus knots in the next section.
Even though in this section we only need $p=1$ specialization of (\ref{Paqt-torus}), this is still quite an intricate expression.
On the other hand, the analysis of constraints arising from the action of various differentials leads to the following expression
\be
\P_n ({\bf 3_1};a,q,t) = \sum_{k=0}^{n-1}  a^{n-1} t^{2k} q^{n(k-1)+1} \frac{(q^{n-1},q^{-1})_k(-a t q^{-1},q)_k}{(q,q)_k} \,, \label{Paqt31}
\ee
and one can verify that this is equal to $p=1$ specialization of (\ref{Paqt-torus}). Explicit values of $\P_n ({\bf 3_1};a,q,t)$ following from (\ref{Paqt31}) for low values of $n$ are given in table~\ref{tab-P31}. Again, note that they are polynomials with positive coefficients, as expected from (\ref{superPdef}).

\begin{table}[h]
\centering
\begin{tabular}{|@{$\Bigm|$}c|@{$\Bigm|$}l|}
\hline
$\, n$ & $\qquad \P_n ({\bf 3_1};a,q,t)$  \\
\hline
\hline
$\, 1$ & $1$   \\
\hline
$\, 2$ & $ a q^{-1} + a q t^2  + a^2 t^3  $    \\
\hline
$\, 3$ & $ a^2 q^{-2} + a^2 q (1 + q) t^2 + a^3 (1 + q) t^3 + a^2 q^4 t^4 +
 a^3 q^3 (1 + q) t^5 + a^4 q^3 t^6$  \\
\hline
$\, 4$ & $  a^3 q^{-3} + a^3 q (1 + q + q^2) t^2 + a^4 (1 + q + q^2) t^3 +
 a^3 q^5 (1 + q + q^2) t^4 + $ \\
& $+ a^4 q^4 (1 + q) (1 + q + q^2) t^5 + a^3 q^4 (a^2 + a^2 q + a^2 q^2 + q^5) t^6 + $ \\
& $+  a^4 q^8 (1 + q + q^2) t^7 + a^5 q^8 (1 + q + q^2) t^8 + a^6 q^9 t^9$  \\
\hline
\end{tabular}
\caption{Colored superpolynomial of the ${\bf 3_1}$ knot for $n=1,2,3,4$.
\label{tab-P31} }
\end{table}

\begin{figure}[ht]
\begin{center}
\includegraphics[width=0.4\textwidth]{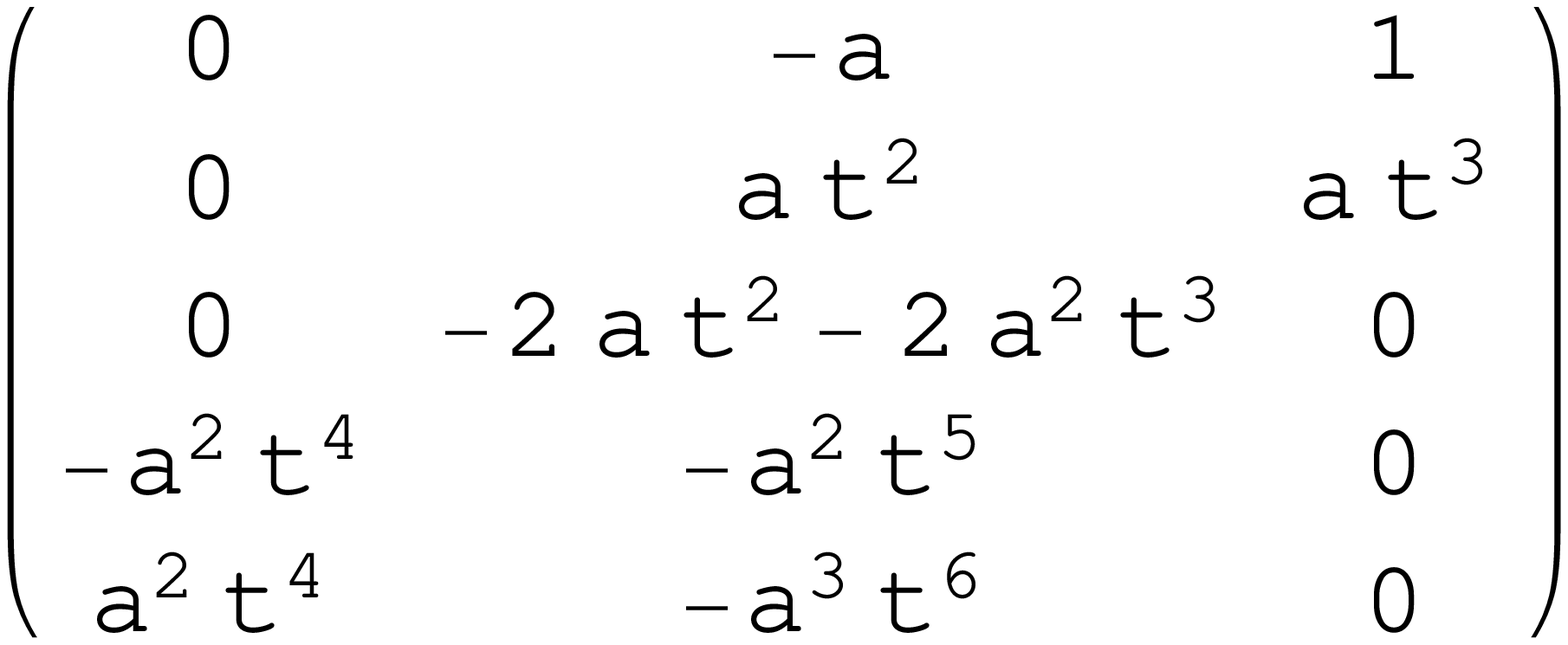}
\caption{Matrix form of the super-$A$-polynomial for the trefoil knot. The conventions are the same as in figure \protect\ref{fig-unknotNewtonMatrix}.}
\label{fig:matrix31}
\end{center}
\end{figure}

\begin{figure}[ht]
\begin{center}
\includegraphics[width=0.4\textwidth]{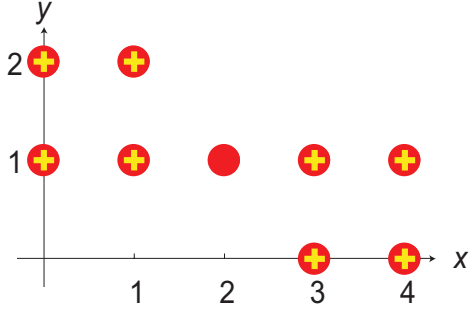}
\caption{Newton polygon of the super-$A$-polynomial for the trefoil knot and its $a=-t=1$ limit. The conventions are the same as in the unknot case in figure \protect\ref{fig-unknotNewtonMatrix}.}
\label{fig:Newton31}
\end{center}
\end{figure}

From the explicit form of the colored superpolynomial (\ref{Paqt31}) we find the recursion
relation it satisfies by using the Mathematica package \textrm{qZeil.m}, see \cite{qZeil}. This recursion relation takes form
\be
\hat A^{\text{super}} (\hat x, \hat y; a,q,t) = a_0 + a_1 \hat{y} + a_2 \hat{y}^2 \,,   \label{Ahat31}
\ee
where
\bea
a_0 & = & \frac{a^2 t^4 (\hat{x}-1) \hat{x}^3 (1 + a q t^3 \hat{x}^2)}{ q (1 + a t^3 \hat{x}) (1 + a t^3 q^{-1} \hat{x}^2) }   \nonumber \\
a_1 & = & -\frac{a (1 + a t^3 \hat{x}^2) \big( q - q^2 t^2 \hat{x} +t^2 (q^2 + q^3 + a t + a q^2 t) \hat{x}^2 + a q^2 t^5 \hat{x}^3 + a^2 q t^6 \hat{x}^4 \big)}{q^2 (1 + a t^3 \hat{x}) (1 + a t^3 q^{-1}\hat{x}^2) }  \nonumber \\
& & \qquad   \nonumber \\
a_2 & = &  1   \nonumber
\eea
The classical super-$A$-polynomial for trefoil knot follows from the $q\to 1$ limit of $\hat A^{\text{super}}$ and reads
\bea
A^{\text{super}} (x, y; a,t) & = & a^2 t^4 (x-1) x^3 + (1 + a t^3 x) y^2 + \label{Asuper31} \\
& & -a\big( 1 - t^2 x + 2 t^2 (1 + a t) x^2 + a t^5 x^3 + a^2 t^6 x^4 \big) y  \,.  \nonumber
\eea
Matrix form of this polynomial is presented in figure \ref{fig:matrix31}, and its Newton polygon is shown in figure \ref{fig:Newton31}.

The same polynomial can be derived from the asymptotic behavior of the colored superpolynomial~(\ref{Paqt31}).
Using integral representation as in (\ref{Pn41-integral}),
\be
\P_n ({\bf 3_1};a,q,t) \; \sim \; \int dz\; e^{\frac{1}{\hbar}\left(\widetilde{\mathcal{W}}({\bf 3_1};z,x)+{\cal O}(\hbar)\right)},
\ee
with the potential
\bea
& & \widetilde{\mathcal{W}}({\bf 3_1};z,x)  =  - \frac{\pi^2}{6} +\big(\log z + \log a \big)\log x + 2(\log t) (\log z)   \label{V31}\\
& & \quad   + \Li_2( x z^{-1}) - \Li_2(x) + \Li_2(-a t) - \Li_2(-a t z) + \Li_2(z) \,,  \nonumber
\eea
and in the limit (\ref{reflimit}) with $z=e^{\hbar k}=const$, we find that equations (\ref{saddle_point41}) and (\ref{yV}) take form
\be
\left\{\begin{array}{l} 1 = \frac{t^2 x (x - z_0) (1 + a t z_0)}{z_0 (z_0-1)}   \\
y =  \frac{a z_0^2 (x-1)}{(x - z_0)}   \end{array} \right.
\ee
Eliminating $z_0$ from these two equations we reproduce the super-$A$-polynomial (\ref{Asuper31}).
We conclude that both Conjecture~1 and Conjecture~2 hold true for the trefoil knot.

We also note that for $t=-1$ and $a=1$ the super-$A$-polynomial  (\ref{Asuper31}) reduces to
\be
A (x, y) = -(x-1)(y-1)(y+x^3) \,,
\label{AAAtref}
\ee
which reproduces the well known $A$-polynomial for the trefoil knot, including the $(y-1)$ factor
associated with abelian flat connections (and the overall immaterial factor $x-1$).
More generally, under a change of variables
\be
Q=a,\quad \beta=x,\quad \alpha=yQ^{-1}\beta^{-6}\frac{1-Q\beta}{1-\beta} \, ,
\ee
and in $t=-1$ limit, (\ref{Asuper31}) reduces (up to an overall factor) to
\be
A(\alpha,\beta,Q) = (1 - Q \beta) + (\beta^3 - \beta^4 + 2 \beta^5 -2 Q \beta^5 - Q \beta^6 + Q^2 \beta^7) \alpha + (-\beta^9 + \beta^{10}) \alpha^2,
\ee
which reproduces the $Q$-deformed or augmentation polynomial for the trefoil knot found in \cite{NgFramed,AVqdef}. Relations between super-$A$-polynomial, $Q$-deformed $A$-polynomial and augmentation polynomial for torus knots are discussed in much more detail in \cite{superA}.

\subsection{$(2,2p+1)$ torus knots}

As the last class of examples we discuss the entire family of $(2,2p+1)$ torus knots, which are also denoted $T^{(2,2p+1)}$.
The $S^r$-colored superpolynomials for this family can be found both from refined Chern-Simons theory, as well as 
from analysis of differentials. The former approach is described in detail in~\cite{FGS} and it leads (after taking into account appropriate $\gamma$-factors and other subtleties) to the expression for
the reduced colored superpolynomial as the ratio of (refined) Chern-Simons partition functions
in $\S^3$ in presence of a given knot and the unknot
\begin{eqnarray}
\P^{S^r} (T^{(2,2p+1)}; a,q,t)&=& (-1)^{pr}\left( \frac{q_1}{q_2} \right)^{\tfrac{pr}{2}}
\frac{Z_{SU(N)}^{\text{ref}}(\S^3,T^{(2,2p+1)}_{\Lambda^r};q_1,q_2)}{Z_{SU(N)}^{\text{ref}}(\S^3,\unknot_{\Lambda^r};q_1,q_2)}  \label{Paqt-torus} \\
& = & \sum_{\ell=0}^r\frac{(qt^2;q)_{\ell}(-at^3;q)_{r+\ell}(-aq^{-1}t;q)_{r-\ell}(q;q)_{r}}
{(q;q)_{\ell}(q^2t^2;q)_{r+\ell}(q;q)_{r-\ell}(-at^3;q)_{r}}\frac{(1-q^{2\ell+1}t^2)}{(1-qt^2)}
\nonumber \\
&&\quad
\times (-1)^{n-1}a^{-\frac{r}{2}}q^{\frac{3(n-1)}{2}-\ell}t^{-(n-1)p-\ell+\frac{r}{2}}
\left[
(-1)^{\ell}a^{\frac{r}{2}}q^{\frac{r^2-\ell(\ell+1)}{2}}t^{\frac{3r}{2}-\ell}
\right]^{2p+1} \,.    \nonumber
\end{eqnarray}

The second approach, based on analysis of differentials, has been employed in \cite{FGSS} and results in the following form of the superpolynomial
\begin{eqnarray}\label{fort2k}
\!\!\!\!\!\!\!\!{\P}^{S^r}(T^{2,2p+1};a,q,t) &=& a^{pr} q^{-pr}  \sum_{0\le k_p \le \ldots \le k_2 \le k_1 \le r}
\left[\!\begin{array}{c} r\\k_1 \end{array}\!\right]\left[\!\begin{array}{c} k_1\\k_2 \end{array}\!\right]\cdots\left[\!\begin{array}{c} k_{p-1}\\k_p \end{array}\!\right]  \times\\
& & \nonumber \!\!\!\!\!\!\!\!\!\!\!\!\!\!\!\!\!\!\!\! \times \,\,\, q^{(2r+1)(k_1+k_2+\ldots+k_p)-\sum_{i=1}^p k_{i-1}k_i} t^{2(k_1+k_2+\ldots+k_p)} \prod_{i=1}^{k_1}(1+aq^{i-2}t),
\end{eqnarray}
where $k_0=r$. One can check that (\ref{Paqt-torus}) and (\ref{fort2k}) agree up to relatively high values of $r$; it would be nice to find an analytic proof valid for all $r$. We are however convinced that both expressions are equal, and various consequences of this fact -- in particular dualities between different UV descriptions of the corresponding $\mathcal{N}=2$, 3d SUSY theory (see section \ref{sec:phys}) -- were presented in \cite{FGSS}. 

For the purpose of this presentation we will focus on the expression the colored superpolynomial in the form (\ref{Paqt-torus}). 
In the asymptotic limit $\hbar\to 0$ limit we find
\begin{eqnarray}
\P^{S^r} (T^{(2,2p+1)}; a,q,t)
\; \sim \; \int dz\; e^{\frac{1}{\hbar}\left(\widetilde{\mathcal{W}} (T^{(2,2p+1)};z,x)+{\cal O}(\hbar)\right)},
\end{eqnarray}
with the potential
\begin{eqnarray}
&& \widetilde{\mathcal{W}} (T^{(2,2p+1)};z,x) =
 p \log (a) \cdot \log x
-p\log (-t)\cdot \log x+(p+1)\pi i\log x+\log (x^{\frac{1}{2}}z^{-1})\cdot\log t
\nonumber \\
&&\quad\quad\quad\quad\quad\quad\quad
+(2p+1)\Biggl(
\pi i\log z+\frac{1}{2}\left((\log x)^2-(\log z)^2\right)
+\log (x^{\frac{3}{2}}z^{-1})\cdot \log t\Biggr)
\nonumber \\
&&\quad\quad\quad\quad\quad\quad\quad
+{\rm Li}_2(z)-{\rm Li}_2(x)-{\rm Li}_2(t^2z)+{\rm Li}_2(-at^3x)
+{\rm Li}_2(t^2xz)
\nonumber \\
&&\quad\quad\quad\quad\quad\quad\quad
-{\rm Li}_2(-at^3xz)
+{\rm Li}_2(xz^{-1})
-{\rm Li}_2(-atxz^{-1})
+{\rm Li}_2(-at)-{\rm Li}_2(1),
\end{eqnarray}
where $z = q^{\ell}$.

For the above potential $\widetilde{\mathcal{W}} (T^{(2,2p+1)};z,x)$,
the critical point condition can simply be expressed as $1=\exp\left(z\partial \widetilde{\mathcal{W}}/\partial z\right)|_{z=z_0}$:
\begin{eqnarray}
1 \; = \; -\frac{t^{-2-2p}(x-z_0)z_0^{-1-2p}(-1+t^2z_0)(1+ at^3 xz_0)}{(-1+z_0)(atx+z_0)(-1 +  t^2 x z_0)} \,.
\label{braid_saddle1}
\end{eqnarray}
and
\begin{eqnarray}
y(x,t,a)&=&\exp\left(x\frac{\partial \widetilde{\mathcal{W}} (T^{(2,2p+1)};z_0,x)}{\partial x}\right)
\nonumber \\
&=&\frac{a^p t^{2 + 2 p} (-1 + x) x^{1 + 2 p} (atx + z_0) (1 + a t^3 x z_0)}{(1 + a t^3 x) (x - z_0) (-1 + t^2 x z_0)} \,.
\label{braid_saddle2}
\end{eqnarray}
Eliminating $z_0$ from the above equations, we find the
super-$A$-polynomial $A^{\text{super}} (x,y;a,t)$ for any $(2,2p+1)$ torus knot.
For small values of $p$, the resulting super-$A$-polynomials are listed
in table \ref{table_super} (where we omitted the extra
factors which appear in the elimination, and picked up the factor that
includes the non-abelian branch of the $SL(2)$ character variety).
In particular, for $p=1$, we obtain (up to an inessential overall factor) the same super-$A$-polynomial for the trefoil (\ref{Asuper31}), which was derived in the previous section starting from another expression for colored superpolynomials (\ref{Paqt31}).
\begin{table}[h]
\centering
\begin{tabular}{|@{$\Bigm|$}c|@{$\Bigm|$}l|}
\hline
\textbf{Knot}  & $A_K^{\text{super}} (x,y;a,t)$  \\
\hline
\hline
$T^{(2,3)}$ & {\scriptsize $y^2 +$}$\frac{1}{1+ a t^3 x}${\scriptsize
     $a(-1 + t^2 x -2 t^2 x^2 - 2 a t^3 x^2 - a t^5 x^3 - a^2 t^6 x^4)y
 +$}$\frac{(x-1)a^2 t^4 x^3}{1 + a t^3 x}$    \\
\hline
$T^{(2,5)}$ & {\scriptsize $y^3
-$}$\frac{a^2}{1 + a t^3 x}$ {\scriptsize $
(1 - t^2 x + 2 t^2 x^2 + 2 a t^3 x^2 - 2 t^4 x^3 - 2 a t^5 x^3 +
 3 t^4 x^4 + 4 a t^5 x^4 + a^2 t^6 x^4 + a t^7 x^5
$
}
\\
&{\scriptsize
$\quad\quad
- a^2 t^8 x^5 +
 2 a^2 t^8 x^6)y^2$
}
\\
&{\scriptsize $+$}$\frac{a^4t^6 (-1 + x) x^5}{(1 + a t^3 x)^2}$
{\scriptsize $(
2 - t^2 x + a t^3 x + 3 t^2 x^2 + 4 a t^3 x^2 + a^2 t^4 x^2 +
 2 a t^5 x^3 + 2 a^2 t^6 x^3 + 2 a^2 t^6 x^4
$
}
\\
&{\scriptsize
$\quad
+ 2 a^3 t^7 x^4 +
 a^3 t^9 x^5 + a^4 t^{10} x^6
)y$
}
\\
&{\scriptsize $-$}$\frac{a^6t^{12} (-1 + x)^2 x^{10}}{(1 + a t^3 x)^2}$  \\
\hline
$T^{(2,7)}$ &
{\scriptsize $y^4
-$}$
\frac{a^3}{1+at^3x}$
{\scriptsize
$(1 - t^2 x + 2 t^2 x^2 + 2 a t^3 x^2 - 2 t^4 x^3 - 2 a t^5 x^3 +
  3 t^4 x^4 + 4 a t^5 x^4 + a^2 t^6 x^4 - 3 t^6 x^5$}
\\&{\scriptsize $\quad\quad
- 4 a t^7 x^5 -
  a^2 t^8 x^5 + 4 t^6 x^6 + 6 a t^7 x^6 + 2 a^2 t^8 x^6 + a t^9 x^7 -
  2 a^2 t^{10} x^7 + 3 a^2 t^{10} x^8
)y^3
$}
\\&{\scriptsize $+$}$\frac{a^6t^8 (-1 + x) x^7}{(1 +a t^3 x)^2}$
{\scriptsize $(
3 - 2 t^2 x + a t^3 x + 6 t^2 x^2 + 8 a t^3 x^2 + 2 a^2 t^4 x^2 -
 3 t^4 x^3 - 2 a t^5 x^3 + a^2 t^6 x^3 + 6 t^4 x^4
$}\\
&{\scriptsize $\quad
+ 12 a t^5 x^4 +
 10 a^2 t^6 x^4 + 4 a^3 t^7 x^4 + 3 a t^7 x^5 + 2 a^2 t^8 x^5 -
 a^3 t^9 x^5 + 6 a^2 t^8 x^6 + 8 a^3 t^9 x^6 + 2 a^4 t^{10} x^6
$}\\
&{\scriptsize $\quad
+  2 a^3 t^{11} x^7 - a^4 t^{12} x^7 + 3 a^4 t^{12} x^8)y^2
$}
\\
&{\scriptsize $-$}$\frac{a^9t^{16} (-1 + x)^2 x^{14}}{(1 +a t^3 x)^3}$
{\scriptsize
$
(3 - t^2 x + 2 a t^3 x + 4 t^2 x^2 + 6 a t^3 x^2 + 2 a^2 t^4 x^2 +
  3 a t^5 x^3 + 4 a^2 t^6 x^3 + a^3 t^7 x^3
$}\\
&{\scriptsize $\quad
+ 3 a^2 t^6 x^4 +
  4 a^3 t^7 x^4 + a^4 t^8 x^4 + 2 a^3 t^9 x^5 + 2 a^4 t^{10} x^5 +
  2 a^4 t^{10} x^6 + 2 a^5 t^{11} x^6 + a^5 t^{13} x^7 + a^6 t^{14} x^8)y$
}
\\
&{+$\frac{a^{12}t^{24} (-1 + x)^3 x^{21}}{(1 + t^3 x)^3}$}  \\
\hline
\end{tabular}
\caption{
Super-$A$-polynomials for $(2,2p+1)$ torus knots with $p=1,2,3$. \label{table_super} }
\end{table}

For $a=1$ the above super-$A$-polynomials reduce to the refined $A$-polynomials of \cite{FGS}. On the other hand,
for $t=-1$ we find the $Q$-deformed $A$-polynomial of \cite{AVqdef}
if the following identification of parameters is performed
\begin{eqnarray}
Q=a,\quad \beta=x,\quad \alpha=yQ^{-p}\beta^{-4p-2}\frac{1-Q\beta}{1-\beta} \, ,
\label{super2av_torus}
\end{eqnarray}
and a related transformation reveals the form of the augmentation polynomial of \cite{NgFramed}.
Precise derivation of the above variable change, as well as explicit relations between super-$A$-polynomial, the augmentation polynomial and $Q$-deformed $A$-polynomial, are discussed in detail in \cite{superA}.


\section{Quantizability}
\label{sec:quantizability}

In this section we discuss the super-$A$-polynomials that we found from the viewpoint of quantizability,
by which we mean the following. For the Conjecture~1 to be formulated in a consistent way,
we must ensure that the leading term $\int \log y \frac{dx}{x}$ in the integral (\ref{VCsuper}) is well-defined,
{\it i.e.} does not depend on the choice of the integration path on the algebraic curve~(\ref{supercurve}).
As explained in~\cite{Apol,abmodel},
this requirement imposes the following constraints on the periods of the imaginary and real parts of $\log y \frac{dx}{x}$, respectively,
\bea
\oint_{\gamma} \Big( \log |x| d ({\rm arg} \, y) - \log |y| d ({\rm arg} \, x) \Big) & = & 0 \,, \label{qcond0} \\
\frac{1}{4 \pi^2} \oint_{\gamma} \Big( \log |x| d \log |y| + ({\rm arg} \, y) d ({\rm arg} \, x) \Big) & \in & \mathbb{Q} \, ,
\label{qcondQ}
\eea
for all {\it closed} paths $\gamma$ on the curve (\ref{supercurve}).
It turns out that these conditions can be further reformulated and interpreted in a variety of ways.
On one hand, it is amusing to observe that the integrand $\eta(x,y)=\log |x| d ({\rm arg} \, y) - \log |y| d ({\rm arg} \, x)$
in \eqref{qcond0} is the image of the symbol $\{ x,y \} \in K_2 (\cC)$ under so-called regulator map,
thereby constituting an immediate link to algebraic K-theory  \cite{Beilinson,Bloch,Fernando}.
As discussed in~\cite{abmodel}, from this K-theory viewpoint the condition that the curve is quantizable
can be rephrased simply as the requirement that $\{ x,y \} \in K_2 (\C (\cC))$ is a torsion class.
On the other hand, this more abstract condition also translates to the down-to-earth statement
that quantizability of the curve requires its defining polynomial to be {\it tempered}.

By definition, a polynomial $A(x,y)$ is tempered if all roots of all face polynomials of its Newton polygon are roots of unity.
Face polynomials are constructed as follows: we need to construct a Newton polygon
corresponding to $A(x, y) = \sum_{i,j} a_{i,j}x^i y^j$, and to each point $(i, j)$ of this polygon we associate the coefficient $a_{i,j}$.
We label consecutive points along each face of the polygon by integers $k=1,2,\ldots$ and,
for a given face, rename monomial coefficients associated to these points as $a_k$.
Then, the face polynomial associated to a given face is defined to be $f(z)=\sum_k a_k z^k$.
Therefore, the quantizability condition requires that all roots of $f(z)$ constructed
for all faces of the Newton polygon must be roots of unity. In what follows we are going
to examine super-$A$-polynomials which we found in examples in section \ref{sec-case} from this perspective.

Ordinary $A$-polynomials have numerical, integer coefficients \cite{CCGL}, and therefore the above
quantization condition imposes certain constraints on values of these coefficients.
For example, the ordinary $A$-polynomial of the figure-eight knot given in (\ref{Afigure8}) satisfies these constraints,
while its close cousin with only slightly different coefficients, discussed {\it e.g.} in \cite{Apol,abmodel}, does not.
Meeting these tight constraints might seem much less trivial in the case of $t$- or $a$-deformed curve,
when coefficients of the defining polynomial depend on these extra parameters.
Nonetheless, this is indeed possible
and the outcome is very simple: the quantization condition implies that both $a$ and $t$ must be roots of unity.
Therefore, even though such $a$ and $t$ cannot be completely arbitrary, they still take values in a dense set of points (on a unit circle).
Below we verify that indeed all super-$A$-polynomials discussed in section~\ref{sec-case}
are tempered (and therefore quantizable) as long as both $a$ and $t$ are roots of unity.
This condition very nicely fits with the fact that in specialization from colored superpolynomial
or HOMFLY polynomial to $sl(N)$ quantum group invariant we substitute $a=q^N$
and in Chern-Simons theory with $SU(N)$ gauge group $q$ is required to be a root of unity,
so that $a=q^N$ is automatically a root of unity as well.

\begin{table}[ht]
\centering
\begin{tabular}{| c || c | }
\hline
\rule{0pt}{5mm}
face & face polynomial \\[3pt]
\hline
\rule{0pt}{5mm}
N  & $ z + a t $ \\[3pt]
NE & $ z + a t^2 $  \\[3pt]
E  & $ a t^2(z+a t^3)^2$    \\[3pt]
SE & $ a^3 t^8(z - a t^2)$  \\[3pt]
S  & $ a^3 t^9 (z + a t) $  \\[3pt]
SW & $ a^2 t^5(z-a t^4) $  \\[3pt]
W  & $ a^2 t^5(z - 1)^2 $  \\[3pt]
NW & $ a t (z-a t^4) $  \\[3pt]
\hline
\end{tabular}
\caption{Face polynomials for the figure-eight knot, corresponding to faces of the octagonal shape formed by non-zero entries of the coefficient matrix in figure \protect\ref{fig:matrix41}. Faces are labeled by compass directions (with N standing for North, \emph{etc}.), with the first row $(0,-at,-1,0)$ of the matrix in figure \protect\ref{fig:matrix41} located in the North.}   \label{tab-face-figure8}
\end{table}

Let us now illustrate the above claim in the examples of various knots discussed in section~\ref{sec-case}.
For each of those knots we construct a Newton polygon and face polynomials of the corresponding super-$A$-polynomials.
In order to construct face polynomials it is convenient to write down a matrix representation of the super-$A$-polynomials.
For instance, for the unknot the Newton polygon and the corresponding matrix representation
are shown in figure~\ref{fig-unknotNewtonMatrix}. In this case, it is clear that roots of face polynomials
are all roots of unity if $a$ and $t$ are roots of unity.
In fact, the unknot is so simple that even a weaker condition is sufficient to hold, namely that the combination $at^3$ is a root of unity.

The matrix coefficients and the Newton polygon for figure-eight knot are given, respectively,
in figures \ref{fig:matrix41} and \ref{fig:Newton41}, and the corresponding face polynomials are presented in table~\ref{tab-face-figure8}.
The face polynomials are manifestly written as products of linear factors, and being tempered requires that both $a$ and $t$ are roots of unity.
An analogous condition holds also for $(2,2p+1)$ torus knots whose Newton polygons have hexagonal shape, and the corresponding face polynomials are given in table \ref{tab-face-torus}. Quantizability conditions are also met for other knots, as verified in \cite{FGSS}. To sum up, from all these examples we conclude that super-$A$-polynomials are quantizable if both $a$ and $t$ are roots of unity; we conjecture that this is the case for all knots. 

\begin{table}[ht]
\centering
\begin{tabular}{| c || c | }
\hline
\rule{0pt}{5mm}
face & face polynomial \\[3pt]
\hline
\rule{0pt}{5mm}
first column & $-(a t^2)^{p(p+1)} (z-1)^p$ \\[3pt]
last column & $(-1)^p (z+a t^3)^p$  \\[3pt]
first row & $z a^p -1$    \\[3pt]
last row & $-(a t^2)^{p(p+1)} \big(z-(a t^2)^{p} \big)$  \\[3pt]
lower diagonal & $(-1)^p  \big(a t^3)^{p} (z - a^{p+1}t^{2p+1} \big)^p$  \\[3pt]
upper diagonal & $(-1)^{p+1} a^p \big(z + a^p t^{2p+2} \big)^p$  \\[3pt]
\hline
\end{tabular}
\caption{Face polynomials for $(2,2p+1)$ torus knots, corresponding to faces of the hexagonal shape formed by non-zero entries of the coefficient matrices for $(2,2p+1)$ torus knots, such as the matrix for the trefoil in figure \protect\ref{fig:matrix41}.}   \label{tab-face-torus}
\end{table}



\section{Interpretation in 3d, $\mathcal{N}=2$ theories}
\label{sec:phys}

The objects we have considered so far, such as super-$A$-polynomials and colored superpolynomials, also have an interesting interpretation in 3d, $\mathcal{N}=2$ SUSY gauge theories. We have already recalled that knot invariants can be described in terms of three-dimensional Chern-Simons theory, and from physics perspective the connection between these two classes of theories arises as a 3d-3d duality associated to complementary compactifications of M5-brane along appropriate three dimensions of its $3+3$ dimensional world-volume \cite{DGH,TYsemiclass,DGG,CCV}. In particular, important properties of both three-dimensional theories (i.e. Chern-Simons and $\mathcal{N}=2$ gauge theory) are encoded in the same algebraic curve. 
From the perspective of $\mathcal{N}=2$ theories this curve plays a role to some extent analogous to the Seiberg-Witten curves of four-dimensional gauge theories \cite{SW-I,SW-II}. 
In what follows we explain that the 3d-3d duality can be extended to incorporate dependence on $a$ and $t$. On Chern-Simons side this corresponds to considering refined Chern-Simons theory with $SL(N)$ gauge group, and on $\mathcal{N}=2$ side these parameters can be interpreted as twisted mass parameters for certain global symmetries $U(1)_{Q}$ and $U(1)_F$. In this context, the algebraic curve mentioned above is precisely the super-$A$-polynomial, and so it carries important information about $\mathcal{N}=2$ theories with those symmetries.

To start with, we recall that the parameter $t$, responsible for
the ``refinement'' or ``categorification'', can be interpreted \cite{FGS} as a twisted mass parameter for the global symmetry $U(1)_F$
in the effective three-dimensional $\cN=2$ theory $T_M$ associated to the knot complement $M = \S^3 \setminus K$:
\be
M \quad \leadsto \quad T_M \,.
\label{TMfromM}
\ee
Moreover, generically, every charged chiral multiplet in a theory $T_M$ contributes to
the effective twisted chiral superpotential a dilogarithm term:
\be
\text{chiral field}~\phi
\qquad \longleftrightarrow \qquad
\begin{array}{l}
\text{twisted superpotential} \\[.1cm]
\Delta \widetilde{\cal W} (\vec x; t) = \textrm{Li}_2
\Big( (-t)^{n_F} \prod_i (x_i)^{n_i} \Big)
\end{array}
\label{xtsuper}
\ee
where $n_F$ is the charge of the chiral multiplet under the global R-symmetry $U(1)_F$
and $\{ n_i \}$ is our (temporary) collective notation for all other charges of $\phi$ under symmetries $U(1)_i$,
some of which may be global flavor symmetries and some of which may be dynamical gauge symmetries,
depending on the problem at hand.\footnote{Below we shall return to the different role of gauge and global symmetries,
but for now we wish to point out a simple rule of thumb that one can read off the matter content
of the theory $T_M$ by counting dilogarithm terms in the function $\widetilde{\cal W} (\vec x; t)$.}
In particular, in the former case, the vev of the corresponding twisted chiral multiplet
is usually called the twisted mass parameter $\tilde m_i = \log x_i$,
of which $\tilde m_F = \log (-t)$ is a prominent example.

\begin{figure}[ht]
\centering
\includegraphics[width=2in]{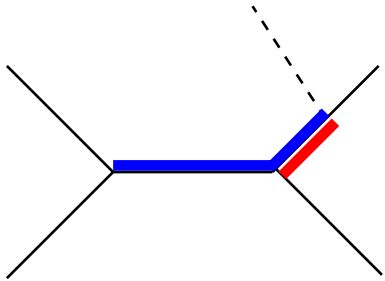}
\caption{A toric Lagrangian brane in the conifold bounds two holomorphic disks
(shown by red and blue intervals in the base of the toric geometry).}
\label{fig:rbbrane}
\end{figure}

The second commutative deformation parameter $a$ also admits a similar interpretation
as a twisted mass parameter for a global symmetry that we denote $U(1)_Q$:
\be
\log a \; = \; \tilde m_Q \,.
\label{mftrel}
\ee
In fact, in the case of the $a$-deformation this interpretation is even more obvious
and can be easily seen in the brane picture, where it corresponds to one of the K\"ahler moduli
of the underlying Calabi-Yau geometry $X$.
For example, the effective low-energy theory on a toric brane in the conifold geometry
has two chiral multiplets that come from two open BPS states shown in blue and red in figure~\ref{fig:rbbrane}.

\bigskip

In this example, the symmetry $U(1)_Q$ responsible for the $a$-deformation comes from
the 2-cycle in the conifold geometry $X$. (The corresponding gauge field $A_{\mu}$
comes from the Kaluza-Klein reduction of the RR 3-form field, $C \sim A \wedge \omega$,
and becomes the starting point for the geometric engineering of $\cN=2$ gauge theories in four dimensions \cite{engineering}.)
In a basis of refined open BPS states shown in figure~\ref{fig:rbbrane},
one state is charged under the symmetry $U(1)_Q$, while the other state is neutral.
Therefore, the effective twisted superpotential $\widetilde{\cal W} (x; a, t)$ of the corresponding
model has two dilogarithm terms, one of which depends on $a$ and the other does not.

Returning to the general theory $T_M$, now we are ready to explain the connection between
the twisted superpotential in this theory and the algebraic curve \eqref{supercurve}
defined as the zero locus of the super-$A$-polynomial.
Roughly speaking, the curve \eqref{supercurve} describes the SUSY vacua in the $\cN=2$ theory $T_M$.
To make this more precise, we need to recall that among the parameters $x_i$ in \eqref{xtsuper}
some correspond to vevs of dynamical fields (and, therefore, need to be integrated out)
and some are twisted masses for global flavor symmetries.
To make the distinction clearer, let us denote the former by $z_i$ (instead of $x_i$),
so that the vevs of dynamical twisted chiral superfields are $\sigma_i = \log z_i$.
Then, in order to find SUSY vacua of the theory $T_M$ we need to extremize $\widetilde{\cal W}$
with respect to these dynamical fields,
\be
\frac{\partial \widetilde{\cal W}}{\partial z_i} \; = \; 0 \,.
\label{wzextr}
\ee
This is exactly what we did {\it e.g.} in \eqref{saddle_point41} when we extremized the potential function \eqref{V41}
for the figure-eight knot ({\it cf.} also \eqref{V31} and \eqref{braid_saddle1} for the case of $(2,2p+1)$ torus knots).
Solving these equations for $z_i$ and substituting the resulting values back into $\widetilde{\cal W}$
gives the effective twisted superpotential, $\widetilde{\cal W}_{\text{eff}}$, that depends only on twisted mass
parameters associated with global symmetries of the $\cN=2$ theory $T_M$.

Besides the symmetries $U(1)_F$ and $U(1)_Q$ which are responsible for $t$- and $a$-deformations, respectively,
our $\cN=2$ theories $T_M$ come with additional global flavor symmetries, one for each component of the link $K$
(or, more generally, one for every torus boundary of $M$).
In particular, if $K$ is a knot --- which is what we assume throughout the present paper --- then,
in addition to $U(1)_F$ and $U(1)_Q$, there is only one extra global symmetry $U(1)_L$
with the corresponding twisted mass parameter that we simply denote $\tilde m$;
it is $x = e^{\tilde m}$ that shortly will be identified with the variable by the same name in the super-$A$-polynomial.
In the brane model,
\be
\begin{matrix}
{\mbox{\rm space-time:}} & \qquad & \R^4 & \times & X \\
& \qquad & \cup &  & \cup \\
{\mbox{\rm D4-brane:}} & \qquad & \R^2 & \times & L
\end{matrix}
\label{surfeng}
\ee
this symmetry $U(1)_L$ can be identified with the gauge symmetry on the D4-brane supported on the Lagrangian
submanifold $L \subset X$. The corresponding gauge field is dynamical when $L$ has finite volume,
while for non-compact $L$ (of infinite volume) the symmetry $U(1)_L$ is a global symmetry.
Moreover, the other global symmetry $U(1)_F$ that plays an important role
in our discussion also can be identified in the brane setup \eqref{surfeng}:
it corresponds to the rotation symmetry of the normal bundle of $\R^2 \subset \R^4$.

To summarize our discussion so far, we can incorporate $U(1)_Q$ and $U(1)_L$ charges in \eqref{xtsuper}
and write the contribution of a chiral multiplet $\phi \in T_M$ to the twisted superpotential as
\be
\text{chiral field}~\phi
\qquad \longleftrightarrow \qquad
\begin{array}{l}
\text{twisted superpotential} \\[.1cm]
\Delta \widetilde{\cal W} (x, z_i; a, t) = \textrm{Li}_2
\Big( a^{n_Q} (-t)^{n_F} x^{n_L} \prod_i (z_i)^{n_i} \Big)
\end{array}
\label{axtsuper}
\ee
Using this dictionary and dilogarithm identities, such as the inversion formula
$\textrm{Li}_2 (x) = - \textrm{Li}_2 \left( \frac{1}{x} \right) - \frac{\pi^2}{6} - \frac{1}{2} [\log (-x)]^2$, from \eqref{V41}
and \eqref{V31} it is easy to read off
the spectrum of the theory $T_M$ for the trefoil knot and for the figure-eight knot:

\begin{table}[htb]
\be
\begin{array}{l@{\;}|@{\;}ccccc@{\;}|@{\;}c}
\multicolumn{7}{c}{\text{trefoil knot}} \\[.1cm]
& \phi_1 & \phi_2 & \phi_3 & \phi_4 & \phi_5 & \text{parameter} \\\hline
U(1)_{\text{gauge}} & -1 & 0 & 0 & -1 & 1 & z \\
U(1)_{F}            & 0 & 0 & 1 & -1 & 0 & -t \\
U(1)_{Q}            & 0 & 0 & 1 & -1 & 0 & a \\
U(1)_{L}            & 1 & -1 & 0 & 0 & 0 & x
\end{array}
\qquad\qquad
\begin{array}{l@{\;}|@{\;}ccccccc}
\multicolumn{8}{c}{\text{figure-eight knot}} \\[.1cm]
& \phi_1 & \phi_2 & \phi_3 & \phi_4 & \phi_5 & \phi_6 & \phi_7 \\\hline
U(1)_{\text{gauge}} & 0 & -1 & 0 & -1 & 0 & -1 & -1 \\
U(1)_{F}            & 0 & 0 & 1 & -1 & 3 & -3 & 0 \\
U(1)_{Q}            & 0 & 0 & 1 & -1 & 1 & -1 & 0 \\
U(1)_{L}            & -1 & 1 & 0 & 0 & 1 & -1 & 0
\end{array}
\notag \ee
\caption{Spectrum of the $\cN=2$ theory $T_M$ for the trefoil and figure-eight knots.}
\label{tab:charges}
\end{table}

\noindent
The terms of lower transcendentality degree, {\it i.e.} products of ordinary logarithms,
also admit a simple interpretation in three-dimensional $\cN=2$ gauge theory $T_M$.
Notice that, in the collective notations $\{ x_i \}$ for global and gauge symmetries $U(1)_i$
used in \eqref{xtsuper}, the dependence of the twisted superpotential $\widetilde{\cal W}$
on $\log x_i$ is always quadratic, see {\it e.g.} \eqref{V41} and \eqref{V31}.
Such terms correspond to supersymmetric Chern-Simons couplings for $U(1)$ gauge (resp. background flavor) fields:
\be
\frac{k_{ij}}{4 \pi} \int A_i \wedge d A_j + \ldots
\qquad \longleftrightarrow \qquad
\begin{array}{l}
\text{twisted superpotential} \\[.1cm]
\Delta \widetilde{\cal W} (\vec x; a,t) = \frac{k_{ij}}{2} \, \log x_i \cdot \log x_j
\end{array}
\label{superCSTM}
\ee
At this point, we should remind the reader that a given $\cN=2$ theory $T_M$ may admit many dual UV descriptions,
with different number of gauge groups and charged matter fields \cite{DGG}.
However, all of these dual descriptions lead to the same space of supersymmetric moduli (twisted mass parameters)
once all dynamical multiplets are integrated out,
{\it i.e.} once the twisted superpotential is extremized \eqref{wzextr} with respect to all $z_i$.

The resulting ``effective'' twisted superpotential $\widetilde{\cal W}_{\text{eff}} (x;a,t)$ depends
only on the twisted mass parameters associated with the global symmetries $U(1)_L$, $U(1)_Q$, and $U(1)_F$.
Then, the algebraic curve \eqref{supercurve} defined as the zero locus of the super-$A$-polynomial
is simply a graph of the function $x \frac{\partial \widetilde{\cal W}_{\text{eff}}}{\partial x}$,
which in a circle compactification of the theory $T_M$ is interpreted as the effective FI parameter:
\be
\boxed{\phantom{\oint} \cM_{\text{SUSY}}: \quad A^{\text{super}} (x,y;a,t) = 0
\qquad \Leftrightarrow \qquad  \log y = x \partial_{x} \widetilde{\cal W}_{\text{eff}} (x; a, t)
\phantom{\oint} }
\label{aviaw}
\ee
We note that an amusing example of dualities between various UV descriptions of the same $\cN=2$ theory, associated to arbitrary $(2,2p+1)$ torus knot, has been analyzed in detail in \cite{FGSS}. In this case two distinct UV descriptions arise from two different of colored superpolynomials given in 
(\ref{Paqt-torus}) and (\ref{fort2k}).






\acknowledgments{We thank Sergei Gukov for a very nice and fruitful collaboration 
that led to the discovery of the super-$A$-polynomial
and other results reviewed in this note.
We also thank Hidetoshi Awata, Satashi Nawata and Marko Sto$\check{\text{s}}$i$\acute{\text{c}}$
for collaborations on these and related topics. 
The work of H.F. is supported by the Grant-in-Aid for Young Scientists
(B) [\# 21740179] from the Japan Ministry of Education, Culture, Sports,
Science and Technology, and the Grant-in-Aid for Nagoya University
Global COE Program, ``Quest for Fundamental Principles in the Universe:
from Particles to the Solar System and the Cosmos.''
The work of P.S. is supported by the Foundation for Polish Science.}



\bibliographystyle{JHEP_TD}
\bibliography{abmodel}

\end{document}